\documentclass[12pt]{article}
\usepackage{amsmath,amssymb}
                       
\usepackage{float} 
\usepackage[english]{babel}
\usepackage{mathtools}
\usepackage{url}
\usepackage[colorlinks,linkcolor=blue,anchorcolor=blue,citecolor=blue,backref=page]{hyperref}

\font\smallit=cmti10

\def\ge{\geqslant}
\def\leq{\leqslant}
\def\geq{\geqslant}
\def\({\left(}
\def\){\right)}
\def\fl#1{\left\lfloor#1\right\rfloor}

\usepackage[norefs,nocites]{refcheck}
\long\def\beginFORGET#1\endFORGET{}

 \def\deg{\mathop{\rm deg}\nolimits}

\def\ord{\mathop{\rm ord}\nolimits}

\let\varphi\varphi
\let\varepsilon\varepsilon

\newtheorem{Thm}{Theorem}[section]

\newtheorem{Cor}[Thm]{Corollary}

\newtheorem{Rem}[Thm]{Remark}

\newtheorem{theorem}{Theorem}

\numberwithin{equation}{section}
\numberwithin{table}{section}
\numberwithin{Thm}{section}

\def\qed{{\hskip0pt\unskip\unskip\nobreak\hfil\penalty50
          \hskip1em\hbox{}\nobreak\hfil
           {$\square$}
          \parfillskip=0pt\finalhyphendemerits=0
          \par}\medskip}

\newenvironment{Proof}
               {{\it Proof.}\ }
               {\qed}

               {\qed}

\newcommand{\BF}{{\mathbb{F}}}

\newcommand{\BZ}{{\mathbb{Z}}}

\newcommand{\CH}{{\mathcal H}}

\newcommand{\CJ}{{\mathcal J}}
\newcommand{\CK}{{\mathcal K}}

\newcommand{\CP}{{\mathcal P}}

\newcommand{\CS}{{\mathcal S}}

\newbox\mybox
\def\arrover#1{\mathrel{
       \setbox\mybox=\hbox spread 1.4em
              {\hfil$\scriptstyle#1$\hfil}
       \vbox{\offinterlineskip\copy\mybox
             \hbox to\wd\mybox{\rightarrowfill}}}}

\def\larrover#1{\mathrel{
       \setbox\mybox=\hbox spread 1.4em
              {\hfil$\scriptstyle#1\vphantom{g}$\hfil}
       \vbox{\offinterlineskip\copy\mybox
             \hbox to\wd\mybox{\leftarrowfill}}}}

\def\ontoover#1{\mathrel{
       \setbox\mybox=\hbox spread 1.4em
              {\hfil$\scriptstyle#1\vphantom{g}$\hfil}
       \vbox{\offinterlineskip\copy\mybox
             \hbox to\wd\mybox{\rightarrowfill\hskip-2.8mm
                               $\rightarrow$}}}}
\def\leftontoover#1{\mathrel{
       \setbox\mybox=\hbox spread 1.4em
              {\hfil$\scriptstyle#1\vphantom{g}$\hfil}
       \vbox{\offinterlineskip\copy\mybox
             \hbox to\wd\mybox{$\leftarrow$\hskip-2.8mm
                               \leftarrowfill}}}}

\begin{document}

\begin{center}
{\bf Lower bounds for periods of Ducci sequences}
\vskip 20pt
{\bf Florian Breuer}\\
{\smallit School of Mathematical and Physical Sciences, University of Newcastle, Newcastle,NSW 2308,  Australia}\\
{\tt florian.breuer@newcastle.edu.au}\\ 
\vskip 10pt
{\bf Igor E. Shparlinski}\\
{\smallit School of Mathematics and Statistics, University of New South Wales, 
Sydney, NSW 2052, Australia}\\
{\tt igor.shparlinski@unsw.edu.au}\\ 
\end{center}
\vskip 30pt

\centerline{\bf Abstract}
\noindent
A Ducci sequence is a sequence of integer $n$-tuples obtained by iterating the map
\[
D : (a_1, a_2, \ldots, a_n) \mapsto \(|a_1-a_2|,|a_2-a_3|,\ldots,|a_n-a_1|\).
\] 
Such a sequence is eventually periodic and we denote by $P(n)$ the maximal period 
of such sequences for given $n$. We prove lower bounds for $P(n)$ by counting certain partitions. 

\pagestyle{myheadings}

\thispagestyle{empty}
\baselineskip=15pt
\vskip 30pt

%
%
%


\section*{\normalsize 1. Introduction}
\addtocounter{section}{1}
\setcounter{equation}{0}

Let $n$ be a positive integer. A Ducci sequence is a sequence of integer $n$-tuples obtained by iterating the map
\[
D : \BZ^n \to \BZ^n;
\]
defined as follows:
\[
 D : (a_1, a_2, \ldots, a_n) \mapsto (|a_1-a_2|,|a_2-a_3|,\ldots,|a_n-a_1|).
\]
There is a long literature on Ducci sequences, see for example~\cite{BLM07,BM08,Bre19,CM37,CST05,Clausing,Ehr90,Lud81,MST06,SB18} and the references therein.

Ducci sequences are eventually periodic, and for each $n$ the largest period is denoted by $P(n)$; it is the period of the sequence starting with $(0,0,\ldots,0,1)$. The sequence $P(1), P(2), \ldots$ is entry A038553 in the Online Encyclopedia of Integer Sequences~\cite{OEIS}.
Since $P(2^k)=1$ and $P(2^km)=2^kP(m)$ if $m$ is not a power of 2, by~\cite[Theorem~4]{Ehr90}, we restrict our attention to odd $n$.

The following upper bounds on $P(n)$ are known. 
Denote by $t=\ord_{(\BZ/n\BZ)^*}(2)$ the multiplicative order of $2$ modulo $n$. If there exists an integer $M$
 for which $2^M \equiv -1 \bmod n$, then we say  {\it $n$ is `with $-1$'\/}.  
 The first of the following upper bounds is proved in~\cite{Lud81}, the second in~\cite{Ehr90} and the third in~\cite{Bre19}.
 
 It is convenient to introduce the following quantities
 \begin{equation}
\label{eq:B1 B2}
B_1(n) = 2^t-1 \qquad \text{and}\qquad B_2(n)=n(2^{t/2}-1). 
 \end{equation}

\begin{theorem} 
\label{thm:Pn div B}
Let $n$ be an odd integer, and $t$ the multiplicative order of $2$ modulo $n$. Then,
\begin{enumerate}
	\item $P(n)$ divides $B_1(n) $.
	\item Suppose $n$ is with $-1$, then $P(n)$ divides $B_2(n)$.
	\item Suppose that $n=p^k$ with $p \equiv 5 \bmod 8$ prime and $2$ is a primitive root modulo $p^k$. If the equation $x^2-py^2=-4$ has no solutions in odd integers $x,y\in\BZ$, then $P(n)$ divides $\frac{1}{3}B_2(n)$.
\end{enumerate}
\end{theorem} 

As for lower bounds, the first of the following results is found again in~\cite{Ehr90}, and the remaining ones in~\cite{GS95}.

\begin{theorem} \label{thm:n div Pn}
Let $n$ be an odd integer. Then
\begin{enumerate}
	\item $n$ divides $P(n)$.
	\item $P(n)=n$ if and only if $n=2^r-1$ for some positive integer $r$.
	\item If $n$ is with $-1$, then $P(n) \geq n(n-2)$.
	\item If $n$ is with $-1$, then $P(n) = n(n-2)$ if and only if $n=2^r+1$ for some positive integer $r$.
\end{enumerate}
\end{theorem} 

The goal of the present paper is to prove new asymptotic lower bounds for $P(n)$ in terms of $t$ and $n$. Our starting point is the fact from~\cite{BLM07} that $P(n)$ is the lowest common multiple of multiplicative orders of elements $\zeta+1$, where $\zeta$ is a primitive $n$th root of unity in the finite field~$\BF_{2^t}$.

Since our results require that at least $t> \sqrt{2n}$ holds, in Section~4
we also give a short survey of known results about the size of $t$. 

\vskip 30pt

\section*{\normalsize 2. Multiplicative orders and partitions}
\addtocounter{section}{1}
\setcounter{Thm}{0}
\setcounter{equation}{0}

Let $1\leq a < n$ be an integer prime to $n$.

Consider the set of representatives, chosen in   the interval $[1,n]$,  of the coset $a\langle 2 \rangle \subseteq (\BZ/n\BZ)^*$  
of the multiplicative group $\langle 2 \rangle$ generated by $2$ in the residues ring modulo $n$. That is,
\begin{align*}
\CS_{a,n} := \{j\in\BZ_{>0}:~1\leq j \leq n,  & \ \gcd(j,n)=1, \\
& \;\exists e_j\in\BZ_{\geq 0}, \;j \equiv a2^{e_j} \bmod n \}
\end{align*} 
Its cardinality is $\#\CS_{a,n}=t$.

Next, we consider the set of partitions of numbers $\leq t-1$ into distinct parts from $\CS_{a,n}$:
 \begin{equation}
\label{eq: set P} 
\CP_{a,n} := \{ (u_j)_{j\in \CS_{a,n}} \in \{0,1\}^{t} \; | \; \sum_{j\in \CS_{a,n}}u_j j \leq t-1 \}. 
 \end{equation}

Our main result is
\begin{Thm}\label{main} 
Suppose $n$ is odd and $a$ is relatively prime to $n$. Then
$P(n) \geq \#\CP_{a,n}$.
\end{Thm}

\begin{Proof}
It follows from~\cite[Theorem~3.9]{BLM07} that $P(n)$ is the lowest common multiple of the multiplicative orders of $\zeta +1$, where $\zeta$ ranges over all $n$th roots of unity $1\neq \zeta\in\BF_{2^t}$.

Let $\zeta \in \BF_{2^t}$ be a primitive $n$th root of unity. 
The idea is to show that every partition in $\CP_{a,n}$ leads to a distinct power of $\zeta + 1$. For this we follow the strategy of~\cite{ASV10}.

Let $u=(u_j)_{j\in \CS_{a,n}}\in \CP_{a,n}$, and set 
\[
Q_u = \sum_{j\in \CS_{a,n}}u_j2^{e_j},
\] 
where $j\equiv a2^{e_j} \bmod n$. We also choose an integer $b$ for which $ab\equiv 1 \bmod n$. Now
\begin{align*}
(\zeta + 1)^{Q_u}  
& =  \prod_{j\in \CS_{a,n}}(\zeta + 1)^{u_j2^{e_j}} 
 =  \prod_{j\in \CS_{a,n}}(\zeta^{2^{e_j}} + 1)^{u_j} \\
& =  \prod_{j\in \CS_{a,n}}(\zeta^{bj} + 1)^{u_j} 
  =   \prod_{j\in \CS_{a,n}}(\vartheta^j + 1)^{u_j},  
\end{align*}
where $\vartheta = \zeta^b \in \BF_{2^t}$ is another primitive $n$th root of unity.

Let 
\[
v=(v_j)_{j\in \CS_{a,n}}\in \CP_{a,n}
\] 
be another partition distinct from $u$, we must show that $v$ gives rise to a distinct power of $\zeta + 1$. 
Suppose $(\zeta+1)^{Q_u} = (\zeta+1)^{Q_v}$, so
\[
\prod_{j\in \CS_{a,n}}(\vartheta^j + 1)^{u_j} = \prod_{j\in \CS_{a,n}}(\vartheta^j + 1)^{v_j}.
\]

Denote by $f(X)\in\BF_2[X]$ the minimal polynomial of $\vartheta$; it has degree $t$. Then $f(X)$ must divide 
$U(X) - V(X)$, where
\[
U(X) = \prod_{j\in \CS_{a,n}}(X^j + 1)^{u_j} \quad\text{and}\quad 
V(X) = \prod_{j\in \CS_{a,n}}(X^j + 1)^{v_j}.
\]
Since these polynomials have degree $\leq t-1 < \deg f$ it follows that $U(X)=V(X)$. After removing common factors from both polynomials (corresponding to $u_j=v_j$), we obtain the identity
 \begin{equation}
\label{eq:prod=prod}
\prod_{h\in \CH}(X^h+1)^{u_h} = \prod_{k\in \CK}(X^k+1)^{v_k}, 
 \end{equation}
where $\CH$ and $\CK$ are disjoint subsets of $\CS_{a,n}$. But now we find that the term of smallest positive degree is $x^e$ where $e$ is the smallest element of $\CH\cup \CK$, but this only appears on one side of the identity~\eqref{eq:prod=prod}. This    contradiction  concludes the proof. 
\end{Proof}

\begin{Rem} Some parts of the proof of  Theorem~\ref{main} can be shortened by appealing  to~\cite[Lemma~1]{Pop14}, 
however for completeness and since~\cite{Pop14} may not be easily accessible, we present a full self-contained proof. 
\end{Rem}

\vskip 30pt

\section*{\normalsize 3. Counting partitions} 
\addtocounter{section}{1}
\setcounter{Thm}{0}
\setcounter{equation}{0}

Now we construct lower bounds for the cardinality of $\CP_{a,n}$ for $n$ 
of prescribed arithmetic structure.
As we have mentioned, these bounds are only useful if $t$ is not too small, specifically $t > \sqrt{2n}$.

Suppose first that $t=\varphi(n)$, that is,  $2$ is a primitive root modulo $n$. In this case, $n=p^k$ must be a power of an odd prime $p$. 

When $n=p$, we find that $\CP_{a,n}$ contains the set of partitions of $n-2$ into distinct parts, and the standard asymptotic for that gives (see e.g. \cite[Theorem 6.4]{And76})

\begin{Cor}\label{Cor1}
Suppose $n=p$ is an odd prime and $2$ is a primitive root modulo $p$. Then, as $n\to\infty$,
\begin{align*}
P(n)  & \geq \exp\left[\left(\frac{\pi}{\sqrt{3}}+o(1)\right)\sqrt{n}\right].
\end{align*}
\end{Cor}

The  case of Corollary~\ref{Cor1} is already contained in~\cite[Theorem~1]{Pop12}; in particular, the completely explicit lower bound (for $2$ a primitive root modulo $n=p$)
\[
P(n) \geq \big(80(n-2)\big)^{-\sqrt{2}} \exp\left( \pi\sqrt{\frac{n-2}{3}}\right)
\]
follows from~\cite[Corollary~4]{Pop12}, see also~\cite{Pop14} for some related results.

Next, suppose that $n=p^k$ and $2$ is a primitive root modulo $n$. For this it suffices that $2$ is a primitive root modulo $p$ and $p$ is not a Wieferich prime, that is,  $2^{p-1} \not\equiv 1 \bmod p^2$.

We have $t=p^{k-1}(p-1)$ and $\CP_{a,n}$ contains the set of partitions of $t-1$ into distinct parts which are not divisible by $p$. An asymptotic formula for the number of such partitions appears in~\cite[Corollary~7.2]{Hag64}, and we obtain

\begin{Cor}\label{Cor2}
Fix an odd non-Wieferich prime $p$ and suppose that $2$ is a primitive root modulo~$p$.
Let $n=p^k$, then as $k\to\infty$, we have
\begin{align*}
P(n)  & \geq   \exp\left[\left(\frac{\pi}{\sqrt{3}}\sqrt{\frac{p-1}{p}}+o(1)\right)\sqrt{n}\right].
\end{align*}
\end{Cor}

If $t<\varphi(n)$, then, inspired by~\cite{GS98}, we estimate the cardinality of $\CP_{a,n}$ as follows. Let $2\leq N < t$ be an integer, and denote by $\CS_{a,n}(N) = \CS_{a,n}\cap[1,N]$. Each subset $\CJ \subseteq \CS_{a,n}(N)$ of cardinality $\#\CJ = J \leq t/N$ produces a valid partition $u\in \CP_{a,n}$, where $u_j = 1$ if $j\in \CJ$ and $u_j=0$ otherwise. Thus we obtain
\[
\#\CP_{a,n} \geq \sum_{J \leq t/N} \left(\begin{array}{c} \#\CS_{a,n}(N) \\ J\end{array}\right).
\]
It remains to estimate $\#\CS_{a,n}(N)$ and choose suitable $a$ and $N$.

It is well known that,  
\[
\#\{j \; : \; 1 \leq j \leq N, \; \gcd(j,n)=1\} = N\varphi(n)/n + O(n^{o(1)}),
\]
see, for example,~\cite[Lemma~2.1]{Shp18}.

Now among the cosets of $\langle 2 \rangle \subseteq (\BZ/n\BZ)^*$, at least one must have at least the average number of representatives in $[1,N]$, so there exists an integer $a$, prime to $n$, for which 
\begin{align*}
\#\CS_{a,n}(N) &\geq \frac{t}{\varphi(n)}\cdot \#\{j \; : \; 1 \leq j \leq N, \; \gcd(j,n)=1\}\\
&=  \frac{t}{\varphi(n)}\( N\varphi(n)/n + O(n^{o(1)})\) =
\(1 + o(1)\) \frac{tN}{n}  
\end{align*} 
as $n\to \infty$, provided $N \ge n^\varepsilon$ for some fixed $\varepsilon > 0$.

Now we choose $N=\fl{\sqrt{2n}}$.  
Since $t \ge n^{1/2 + \varepsilon}$, we have 
\[
\#\CS_{a,n}(N) \ge \frac{tN}{n} + O(n^{o(1)}) = \(2+o(1)\) \frac{t}{N}.
\]
Thus by the Stirling formula
\begin{align*}
\#\CP_{a,n}&  \geq \sum_{J \leq t/N} \binom{\#\CS_{a,n}(N)}{J}   
\ge  \binom{\#\CS_{a,n}(N)}{\fl{t/N} } \\
& \ge  \exp\( \(2 \log 2+o(1)\) \frac{t}{N}\).
\end{align*}

Thus we have proved

\begin{Cor}
\label{cor: t vs n} 
Suppose $n$ is odd and $t$ is the multiplicative order of 2 modulo $n$. Then
\[
P(n) \geq \exp\left[\( \log 4 + o(1)\)\frac{t}{\sqrt{2n}}\right].
\]
\end{Cor}

In particular, if $n = p^k$ then it is easy to show that $t \ge c(p) p^k$, 
where $c(p) > 0$ depends only on $p$,
hence Corollary~\ref{cor: t vs n}  
gives a version of  Corollary~\ref{Cor2} in the form
\[
P(n) \geq \exp\( c(p) \sqrt{n}\).
\]

We remark that the condition $t>\sqrt{2n}$ of Corollary~\ref{cor: t vs n} corresponds to the limits of our 
method. Indeed, there are about $\varphi(n)/t$ distinct cosets $\CS_{a,n}$ and since $\varphi(n) = n^{1+o(1)}$ 
each of them is expected to contain very few elements from the interval $[1,t]$ which are the only suitable elements 
which can be used in the construction of the set $\CP_{a,n}$ given by~\eqref{eq: set P}.

Since 
\[
\frac{\log 4}{\sqrt{2}} \approx 0.98025 \qquad \text{and} \qquad \frac{\pi}{\sqrt{3}} \approx 1.8138,
\] 
in the case of $t\approx n$ we recover a result similar to Corollaries~\ref{Cor1} and~\ref{Cor2}, but with a smaller constant in the exponent.

Our lower bounds are quite small compared to the upper bounds $P(n) \leq B_1(n) \sim 2^t$ and $P(n) \leq B_2(n) \sim n 2^{t/2}$, 
see~\eqref{eq:B1 B2}, 
which follow from Theorem~\ref{thm:Pn div B}.  
 On the other hand, they are typically much stronger than linear and quadratic in $n$  
lower bounds which one can extract from Theorem~\ref{thm:n div Pn}.

\vskip 30pt
\section*{\normalsize 4. Lower bounds on multiplicative orders} 
\addtocounter{section}{1}
\setcounter{Thm}{0}
\setcounter{equation}{0}

Since the quality of our bounds depends rather dramatically on the multiplicative order of $2$ modulo $n$, here we give a
short outline of known results. 

First we observe that the applicability of   Corollary~\ref{Cor1} for infinitely many prime $n=p$ is equivalent to 
{\it Artin's conjecture\/}, see~\cite{Mor} for an exhaustive survey.  
On the other hand, we are not aware of any conditional (and certainly unconditional) results or well-established conjectures
towards a version of  Artin's conjecture for non-Wieferich primes which appear Corollary~\ref{Cor2}.
It is natural to expect  that there are infinitely many such primes but known results are scarce~\cite{Silv}. 

Primes $p$ and integers $n$ for which $t$ is large, in particular exceeds $\sqrt{p}$, have been studied 
in many different contexts, but most commonly in the theory of {\it pseudorandom number generators\/}. 
These results originate from the work of Erd{\H o}s and  Murty~\cite{ErdMur} and are conveniently summarised in~\cite{KurPom}. 
For example, for any function $\psi(n) \to 0$ as $n\to \infty$ we have $t\ge n^{1/2+\psi(n)}$ for almost all (in a sense 
of relative density) primes $p=n$ (see~\cite[Theorem~1]{ErdMur}) and odd integers $n$  (see~\cite[Theorem~11]{KurPom}). 
Furthermore, for a positive proportion of primes  $p=n$   (see~\cite[Lemma~19]{KurPom})) and odd integers $n$  (see~\cite[Theorem~21]{KurPom})
 we have $t \ge n^{0.677}$.

\vskip 30pt

\section*{\normalsize 5. Numerical results}
\addtocounter{section}{1}
\setcounter{Thm}{0}
\setcounter{equation}{0}

It is interesting to compare the lower bound of Theorem~\ref{main}  with actual values of $P(n)$. 
Table~\ref{PanValues} 
shows numerical values of $P(n)$ and $\#\CP_{a,n}$ for odd $n\leq 101$ and 
a representative $a$ for each coset $(\BZ/n\BZ)^*/\langle 2 \rangle$.
 Unsurprisingly, the largest value of $\#\CP_{a,n}$ is achieved for $a=1$ in these small cases, due to the presence of small powers of two in $\CS_{a,n}$. 
However, when $n=109$, we find that 
$$\#\CP_{1,109} = 99 < 178 = \#\CP_{3,109}  = \max_{\gcd(a,109)=1} \#\CP_{a,109} .
$$

These values were computed using Sage.
\begin{table}%
\begin{center}
{\tiny
\begin{tabular}{|rrrrr|}
\hline
	$n$  &     $P(n)$  & $t$  & $a$ & $\#\CP_{a,n}$ \\ \hline\hline
   3    &           3 &   2  &  1 &     2 \\ \hline
   5    &          15 &   4  &  1 &     5 \\ \hline
   7    &           7 &   3  &  1 &     3 \\
   -    &           - &   -  &  3 &     1 \\ \hline
   9    &          63 &   6  &  1 &     7 \\ \hline
  11    &         341 &  10  &  1 &    33 \\ \hline
  13    &         819 &  12  &  1 &    55 \\ \hline
  15    &          15 &   4  &  1 &     4 \\
   -    &           - &   -  &  7 &     1 \\ \hline
  17    &         255 &   8  &  1 &     8 \\
   -    &           - &   -  &  3 &     5 \\ \hline
  19    &        9709 &  18  &  1 &   207 \\ \hline
  21    &          63 &   6  &  1 &     6 \\
   -    &           - &   -  &  5 &     2 \\ \hline
  23    &        2047 &  11  &  1 &    28 \\
   -    &           - &   -  &  5 &     4 \\ \hline
  25    &       25575 &  20  &  1 &   190 \\ \hline
  27    &       13797 &  18  &  1 &    79 \\ \hline
  29    &      475107 &  28  &  1 &  1261 \\ \hline
  31    &          31 &   5  &  1 &     5 \\
   -    &           - &   -  &  3 &     2 \\
   -    &           - &   -  &  5 &     1 \\
   -    &           - &   -  &  7 &     1 \\
   -    &           - &   -  & 11 &     1 \\
   -    &           - &   -  & 15 &     1 \\ \hline
  33    &        1023 &  10  &  1 &    10 \\
   -    &           - &   -  &  5 &     3 \\ \hline
  35    &        4095 &  12  &  1 &    16 \\
   -    &           - &   -  &  3 &     4 \\ \hline
  37    &     3233097 &  36  &  1 &  4310 \\ \hline
  39    &        4095 &  12  &  1 &    22 \\
   -    &           - &   -  &  7 &     2 \\ \hline
  41    &       41943 &  20  &  1 &    70 \\
   -    &           - &   -  &  3 &    25 \\ \hline
  43    &        5461 &  14  &  1 &    17 \\
   -    &           - &   -  &  3 &    10 \\
   -    &           - &   -  &  7 &     4 \\ \hline
  45    &        4095 &  12  &  1 &    12 \\
   -    &           - &   -  &  7 &     3 \\ \hline
  47    &     8388607 &  23  &  1 &   241 \\
   -    &           - &   -  &  5 &    14 \\ \hline
  49    &     2097151 &  21  &  1 &    53 \\
   -    &           - &   -  &  3 &    27 \\ \hline
  51    &         255 &   8  &  1 &     8 \\
   -    &           - &   -  &  5 &     3 \\
   -    &           - &   -  & 11 &     1 \\
   -    &           - &   -  & 19 &     1 \\ \hline
  53    &  3556769739 &  52  &  1 & 35680 \\ \hline
  55    &     1048575 &  20  &  1 &    66 \\
   -    &           - &   -  &  3 &     8 \\ \hline
  57    &       29127 &  18  &  1 &    33 \\
   -    &           - &   -  &  5 &     8 \\ \hline
  59    & 31675383749 &  58  &  1 & 72503 \\ \hline
	61   &  65498251203 &  60  &  1 & 91103 \\ \hline
  63   &           63 &   6  &  1 &     6 \\
   -   &            - &   -  &  5 &     2 \\
   -   &            - &   -  & 11 &     1 \\
   -   &            - &   -  & 13 &     1 \\
   -   &            - &   -  & 23 &     1 \\
   -   &            - &   -  & 31 &     1 \\ \hline
\end{tabular}
}
{\tiny
\begin{tabular}{|rrrrr|}
\hline
	$n$  &     $P(n)$  & $t$  & $a$ & $\#\CP_{a,n}$ \\ \hline\hline
  65   &         4095 &  12  &  1 &    12 \\
   -   &            - &   -  &  3 &     4 \\
   -   &            - &   -  &  7 &     3 \\
   -   &            - &   -  & 11 &     2 \\ \hline
  67   & 575525617597 &  66  &  1 & 176945 \\ \hline
  69   &      4194303 &  22  &  1 &    31 \\
   -   &            - &   -  &  5 &    17 \\ \hline
  71   &  34359738367 &  35  &  1 &  1427 \\
   -   &            - &   -  &  7 &    35 \\ \hline
  73   &          511 &   9  &  1 &     9 \\
   -   &            - &   -  &  3 &     3 \\
   -   &            - &   -  &  5 &     3 \\
   -   &            - &   -  &  9 &     1 \\
   -   &            - &   -  & 11 &     1 \\
   -   &            - &   -  & 13 &     1 \\
   -   &            - &   -  & 17 &     1 \\
   -   &            - &   -  & 25 &     1 \\ \hline
  75   &      1048575 &  20  &  1 &    24 \\
   -   &            - &   -  &  7 &     6 \\ \hline
  77   &   1073741823 &  30  &  1 &   100 \\
   -   &            - &   -  &  3 &    70 \\ \hline
  79   & 549755813887 &  39  &  1 &  1028 \\
   -   &            - &   -  &  3 &   106 \\ \hline
  81   &  10871635887 &  54  &  1 &  6159 \\ \hline
  83 & 182518930210733 &   82 &   1 & 911361 \\ \hline
  85   &          255 &   8  &  1 &     8 \\
   -   &            - &   -  &  3 &     3 \\
   -   &            - &   -  &  7 &     2 \\
   -   &            - &   -  &  9 &     1 \\
   -   &            - &   -  & 13 &     1 \\
   -   &            - &   -  & 21 &     1 \\
   -   &            - &   -  & 29 &     1 \\
   -   &            - &   -  & 37 &     1 \\ \hline
  87   &    268435455 &  28  &  1 &   154 \\
   -   &            - &   -  &  5 &     9 \\ \hline
  89   &         2047 &  11  &  1 &    11 \\
   -   &            - &   -  &  3 &     6 \\
   -   &            - &   -  &  5 &     3 \\ 
   -   &            - &   -  &  9 &     2 \\
   -   &            - &   -  & 11 &     1 \\
   -   &            - &   -  & 13 &     1 \\
   -   &            - &   -  & 19 &     1 \\
   -   &            - &   -  & 33 &     1 \\ \hline
  91   &         4095 &  12  &  1 &    12 \\
   -   &            - &   -  &  3 &     8 \\
   -   &            - &   -  &  9 &     2 \\
   -   &            - &   -  & 11 &     2 \\
   -   &            - &   -  & 17 &     1 \\
   -   &            - &   -  & 19 &     1 \\ \hline
  93   &         1023 &  10  &  1 &    10 \\
   -   &            - &   -  &  5 &     2 \\
   -   &            - &   -  &  7 &     2 \\
   -   &            - &   -  & 11 &     1 \\
   -   &            - &   -  & 17 &     1 \\
   -   &            - &   -  & 23 &     1 \\ \hline
  95   &  22906492245 &  36  &  1 &   905 \\
   -   &            - &   -  &  7 &    17 \\ \hline
  97   &   1627389855 &  48  &  1 &  2216 \\
   -   &            - &   -  &  5 &   283 \\ \hline
  99   &      3243933 &  30  &  1 &    49 \\
   -   &            - &   -  &  5 &    32 \\ \hline
	101  & 37905296863701641 & 100 & 1 & 4827382 \\ \hline
\end{tabular}
}
\end{center}
\caption{Values of $P(n)$ and $\#\CP_{a,n}$ for odd $n\leq 101$.}
\label{PanValues}
\end{table}

\vskip 30pt
\newpage
 
\paragraph{Acknowledgments}
The first author thanks the Alexander-von-Humboldt Foundation for support, the Universit\"at Heidelberg for hospitality and Hannes Breuer for interesting discussions.
The second author was  supported in part by the Australian Research Council Grant  DP180100201. 

The authors are grateful to the organisers of the the {\it Sixth  Number Theory Down Under\/}  Conference
(NTDU-6),  Canberra, 24--27 September, 2018, for creating a very encouraging and collaborative atmosphere, 
which has led to this work.



\end{document}